\documentclass{article}

\usepackage{amssymb,amsmath,theorem,euscript}

\newcounter{sec}

\newcounter{punct}[sec]

\def\punct{\refstepcounter{punct}{\arabic{sec}.\arabic{punct}.  }}

\def\COUNTERS{\addtocounter{sec}{1}
              \setcounter{punct}{0}
          \setcounter{equation}{0}
          \setcounter{theorem}{0}
                  }

\newtheorem{theorem}{Theorem}[sec]
\newtheorem{proposition}[theorem]{Proposition}

\newtheorem{lemma}[theorem]{Lemma}

 \def\ov{\overline}
\def\wt{\widetilde}

\begin{document}

\def\OO{\mathrm{O}}
\def\GLO{\mathrm{GLO}}
\def\Coll{\mathrm{Coll}}
\def\kappa{\varkappa}
\def\Mat{\mathrm{Mat}}
\def\U{\mathrm U}

\def\R{\mathbb{R}}
\def\C{\mathbb{C}}

\def\la{\langle}
\def\ra{\rangle}

 \def\cA{\mathcal A}
\def\cB{\mathcal B}
\def\cC{\mathcal C}
\def\cD{\mathcal D}
\def\cE{\mathcal E}
\def\cF{\mathcal F}
\def\cG{\mathcal G}
\def\cH{\mathcal H}
\def\cJ{\mathcal J}
\def\cI{\mathcal I}
\def\cK{\mathcal K}
 \def\cL{\mathcal L}
\def\cM{\mathcal M}
\def\cN{\mathcal N}
 \def\cO{\mathcal O}
\def\cP{\mathcal P}
\def\cQ{\mathcal Q}
\def\cR{\mathcal R}
\def\cS{\mathcal S}
\def\cT{\mathcal T}
\def\cU{\mathcal U}
\def\cV{\mathcal V}
 \def\cW{\mathcal W}
\def\cX{\mathcal X}
 \def\cY{\mathcal Y}
 \def\cZ{\mathcal Z}
%%% END MATHCAL %%%%%%%%%%%%%%%%%%%%%%%%%%%%%%%%% %%%%%%%%%%%%%%%%%%%%%%%%%%%%%%%% %%%
\def\0{{\ov 0}}
% \def\1{{\ov 1}}
 %%%%%%%%%%%%%%%%%%%%%%%%%%%% %%%%%%%%%%%%%%%%%%%%%%%%%%%%%%%%%%% %%% BEGIN GOTIC
 \def\frA{\mathfrak A}
 \def\frB{\mathfrak B}
\def\frC{\mathfrak C}
\def\frD{\mathfrak D}
\def\frE{\mathfrak E}
\def\frF{\mathfrak F}
\def\frG{\mathfrak G}
\def\frH{\mathfrak H}
\def\frI{\mathfrak I}
 \def\frJ{\mathfrak J}
 \def\frK{\mathfrak K}
 \def\frL{\mathfrak L}
\def\frM{\mathfrak M}
 \def\frN{\mathfrak N} \def\frO{\mathfrak O} \def\frP{\mathfrak P} \def\frQ{\mathfrak Q} \def\frR{\mathfrak R}
 \def\frS{\mathfrak S} \def\frT{\mathfrak T} \def\frU{\mathfrak U} \def\frV{\mathfrak V} \def\frW{\mathfrak W}
 \def\frX{\mathfrak X} \def\frY{\mathfrak Y} \def\frZ{\mathfrak Z} \def\fra{\mathfrak a} \def\frb{\mathfrak b}
 \def\frc{\mathfrak c} \def\frd{\mathfrak d} \def\fre{\mathfrak e} \def\frf{\mathfrak f} \def\frg{\mathfrak g}
 \def\frh{\mathfrak h} \def\fri{\mathfrak i} \def\frj{\mathfrak j} \def\frk{\mathfrak k} \def\frl{\mathfrak l}
 \def\frm{\mathfrak m} \def\frn{\mathfrak n} \def\fro{\mathfrak o} \def\frp{\mathfrak p} \def\frq{\mathfrak q}
 \def\frr{\mathfrak r} \def\frs{\mathfrak s} \def\frt{\mathfrak t} \def\fru{\mathfrak u} \def\frv{\mathfrak v}
 \def\frw{\mathfrak w} \def\frx{\mathfrak x} \def\fry{\mathfrak y} \def\frz{\mathfrak z} \def\frsp{\mathfrak{sp}}
 %% This is Lie algebra %%% END GOTIC
%%%%%%%%%%%%%%%%%%%%%%%%%%%%%%%% %%%%%%%%%%%%%%%%%%%%%%%%%%%%%%%%%
%%% BEGIN MATHBF
 \def\bfa{\mathbf a} \def\bfb{\mathbf b} \def\bfc{\mathbf c} \def\bfd{\mathbf d} \def\bfe{\mathbf e} \def\bff{\mathbf f}
 \def\bfg{\mathbf g} \def\bfh{\mathbf h} \def\bfi{\mathbf i} \def\bfj{\mathbf j} \def\bfk{\mathbf k} \def\bfl{\mathbf l}
 \def\bfm{\mathbf m} \def\bfn{\mathbf n} \def\bfo{\mathbf o} \def\bfp{\mathbf p} \def\bfq{\mathbf q} \def\bfr{\mathbf r}
 \def\bfs{\mathbf s} \def\bft{\mathbf t} \def\bfu{\mathbf u} \def\bfv{\mathbf v} \def\bfw{\mathbf w} \def\bfx{\mathbf x}
 \def\bfy{\mathbf y} \def\bfz{\mathbf z} \def\bfA{\mathbf A} \def\bfB{\mathbf B} \def\bfC{\mathbf C} \def\bfD{\mathbf D}
 \def\bfE{\mathbf E} \def\bfF{\mathbf F} \def\bfG{\mathbf G} \def\bfH{\mathbf H} \def\bfI{\mathbf I} \def\bfJ{\mathbf J}
 \def\bfK{\mathbf K} \def\bfL{\mathbf L} \def\bfM{\mathbf M} \def\bfN{\mathbf N} \def\bfO{\mathbf O} \def\bfP{\mathbf P}
 \def\bfQ{\mathbf Q} \def\bfR{\mathbf R} \def\bfS{\mathbf S} \def\bfT{\mathbf T} \def\bfU{\mathbf U} \def\bfV{\mathbf V}
 \def\bfW{\mathbf W} \def\bfX{\mathbf X} \def\bfY{\mathbf Y} \def\bfZ{\mathbf Z} \def\bfw{\mathbf w}
 %%% END MATHBF
%%%%%%%%%%%%%%%%%%%%%%%%%%%%%%% %%%%%%%%%%%%%%%%%%%%%%%%%%%%%%%%%
 %%% BEGIN MATHBB
 \def\R {{\mathbb R }} \def\C {{\mathbb C }} \def\Z{{\mathbb Z}} \def\H{{\mathbb H}}
  \def\K{{\mathbb K}}
   \def\k{{\Bbbk}}
 \def\N{{\mathbb N}} \def\Q{{\mathbb Q}} \def\A{{\mathbb A}} \def\T{\mathbb T} 
 \def\G{\mathbb G}
 \def\bbA{\mathbb A} \def\bbB{\mathbb B} \def\bbD{\mathbb D} \def\bbE{\mathbb E} \def\bbF{\mathbb F} \def\bbG{\mathbb G}
 \def\bbI{\mathbb I} \def\bbJ{\mathbb J} \def\bbL{\mathbb L} \def\bbM{\mathbb M} \def\bbN{\mathbb N} \def\bbO{\mathbb O}
 \def\bbP{\mathbb P} \def\bbQ{\mathbb Q} \def\bbS{\mathbb S} \def\bbT{\mathbb T} \def\bbU{\mathbb U} \def\bbV{\mathbb V}
 \def\bbW{\mathbb W} \def\bbX{\mathbb X} \def\bbY{\mathbb Y} \def\kappa{\varkappa} \def\epsilon{\varepsilon}
 \def\phi{\varphi} \def\le{\leqslant} \def\ge{\geqslant}

\def\P{\mathbf P}

\def\GL{\mathrm {GL}}
\def\bGL{\mathbf {GL}}
\def\GLB{\mathrm {GLB}}

\def\bGr{\mathbf {Gr}}
\def\Gr{\mathrm {Gr}}
\def\bFl{\mathbf {Fl}}

\def\1{\mathbf {1}}

 \newcommand{\Dim}{\mathop {\mathrm {Dim}}\nolimits}
  \newcommand{\codim}{\mathop {\mathrm {codim}}\nolimits}
   \newcommand{\im}{\mathop {\mathrm {im}}\nolimits}
\newcommand{\ind}{\mathop {\mathrm {ind}}\nolimits}
\newcommand{\graph}{\mathop {\mathrm {graph}}\nolimits}

\def\F{\bbF}

\def\sm{\smallskip}

\begin{center}
\Large\bf

On multiplication of double cosets for $\GL(\infty)$ over a finite field

\medskip

\large \sc
Yury A. Neretin%
\footnote{Supported by the grant FWF, P25142.}
\end{center}

{\small
We consider a group $\GL(\infty)$, its parabolic subgroup $B$ corresponding to partition
$\infty=\infty+m+\infty$. Denote by $\P$ the kernel of 
the natural homomorphism $B\to \GL(m)$. We show that the space of double cosets of $\GL(\infty)$
by $\P$ admits a natural structure of a semigroup. In fact this semigroup acts in subspaces of
$\P$-fixed vectors of some unitary representations of $\GL(\infty)$ over finite field.}

\section{Formulation of result}

\COUNTERS

{\bf\punct Notation.}
Let $\k$ be a field.  Denote by $\GL(\infty,\k)$ the group of all infinite matrices $g$ over
$\k$ such that $g-1$ has finite number of non-zero matrix elements. We call such matrices
$g$ {\it finitary}.  We denote:

\sm

--- $\U(\infty)\subset \GL(\infty,\C)$
 the group of finitary unitary matrices over $\C$; 
 
 \sm
 
 --- $\OO(\infty)\subset \GL(\infty,\R)$  the group
 of finitary real orthogonal  matrices;
 
\sm 
 
 ---  $\OO(\infty,\C)\subset \GL(\infty,\C)$
 the group of finitary complex orthogonal matrices.

\sm
 
% By $\F_q$ we denote the finite field with $q$ elements.

Let $G$ be a group, $K$ a subgroup. By $G//K$ we denote the space of conjugacy classes of $G$
with respect to $K$, by $K\setminus G/K$ the double cosets with respect to $K$.

By $\1_k$ we denote the unit matrix of order $k$.

\sm

{\bf\punct Multiplication of conjugacy classes.}  Denote by $\GL(m+\infty,\k)$
the same group $\GL(\infty,\k)$ considered as a group of finitary block matrices
$\begin{pmatrix}A&B\\C&D \end{pmatrix}$
of size $m+\infty$. Denote by $K$ the subgroup, consisting of matrices
of the form
$\begin{pmatrix}\1_m&0\\0&H \end{pmatrix}$, 
$$K\simeq \GL(\infty,\k).$$
 The set
$$G//K=\GL(m+\infty,\C)//\GL(\infty,\C)$$
is a semigroup with respect to the following $\circ$-multiplication.
Consider two matrices 
$\begin{pmatrix}A&B\\C&D \end{pmatrix}$, $\begin{pmatrix}P&Q\\R&T \end{pmatrix}$.
We  split them as block matrices
of size $m+(N+\infty)$. For sufficiently large $N$ we get matrices of the following
block structure:
$$
\begin{pmatrix}A&B\\C&D \end{pmatrix}=\begin{pmatrix}a&|&b&0\\
\text{---}& &\text{---}&\text{---}\\
c&|&d&0\\0&|&0&\1_\infty \end{pmatrix}, \qquad
\begin{pmatrix}P&Q\\R&T \end{pmatrix}=\begin{pmatrix}p&|&q&0\\
\text{---}& &\text{---}&\text{---}\\
r&|&t&0\\0&|&0&\1_\infty \end{pmatrix}
$$
Then
\begin{multline}
\begin{pmatrix}a&b&0\\c&d&0\\0&0&\1_\infty \end{pmatrix}\circ
\begin{pmatrix}p&q&0\\r&t&0\\0&0&\1_\infty \end{pmatrix}
:=
\begin{pmatrix}a&b&0&0\\c&d&0&0\\0&0&\1_N&0\\ 0&0&0&\1_\infty \end{pmatrix}
 \begin{pmatrix}p&0&q&0\\0&\1_N&0&0\\r&0&t&0\\
 0&0&0&\1_\infty \end{pmatrix}
 =\\=
 \begin{pmatrix}ap&|&b&aq&0\\ 
\text{---}& &\text{---}&\text{---}&\text{---}\\
 cp&|&d&cq&0\\r&|&0&t&0
 \\ 0&|&0&0&\1_\infty \end{pmatrix}
 \label{eq:circ}
\end{multline}
The matrix in the right-hand side has size $m+N+N+\infty$. We unite 
$N+N+\infty=\infty$ and get a matrix of size $m+\infty$.

\begin{theorem}
{\rm a)}
The $\circ$-multiplication is a well-defined associative operation on the set
of conjugacy classes
$\GL(m+\infty,\k)//\GL(\infty,\k)$.

\sm

{\rm b)} The $\circ$-multiplication is a well-defined associative operation on the set
of conjugacy classes
$\U(m+\infty)//\U(\infty)$. 
\end{theorem}

The statement (if it is formulated) is more-or-less obvious.
Various versions of these  semigroups are classical topics of system theory and operator theory,
see, e.g., \cite{Bro}, \cite{Dym}, \cite{GGK}, \cite{Haz}, \cite{RR}. 

\sm

{\bf\punct Multiplications of double cosets.}

\begin{theorem}
{\rm a)} The formula {\rm(\ref{eq:circ})} determines an associative operation on double cosets
$\OO(\infty)\setminus \U(m+\infty)/\OO(\infty)$.

\sm

{\rm b)} The formula {\rm(\ref{eq:circ})} determines an associative operation on double cosets
$\OO(\infty,\C)\setminus \GL(m+\infty,\C)/\OO(\infty,\C)$.
\end{theorem}

The  semigroup $\OO(\infty)\setminus \U(m+\infty)/\OO(\infty)$ arises in  representation
theory in the following situation (\cite{Olsh-GB}, \cite{Ner-book}).
Consider a unitary representation $\rho$ of the group
$\U(m+\infty)$. Denote by $H_m$ the subspace of $\OO(\infty)$-fixed  vectors, by $\Pi_m$
the operator of orthogonal projection to $H_m$. Let
 $\gamma\in \OO(\infty)\setminus \U(m+\infty)/\OO(\infty)$ be a double coset, $g\in \U(m+\infty)$
 be a representative of $\gamma$. We define an operator 
 $$
 \ov\rho(\gamma):H_m\to H_m
 $$
by
$$
\ov\rho(\gamma)=\Pi_m \rho(g)\Bigr|_{H_m}
$$
It is easy to show that $\rho(\gamma)$ does not depend on the choice of a representative $g\in\gamma$.

\begin{theorem}
For any $\gamma_1$, $\gamma_2\in \OO(\infty)\setminus \U(m+\infty)/\OO(\infty)$
the following identity {\rm(}'multiplicativity theorem'{\rm)} holds
$$
\rho(\gamma_1)\rho(\gamma_2)=
\rho(\gamma_1\circ \gamma_2)
.
$$
\end{theorem}

See \cite{Olsh-GB}; for details, see \cite{Ner-book},  Section IX.4. According
\cite{Olsh-GB}, theorems of
this type hold for all infinite-dimensional limits $G(\infty)\supset K(\infty)$
of symmetric pairs $G(n)\supset K(n)$.
Recently \cite{Ner-char}, \cite{Ner-faa}, \cite{Ner-invariant},
it was observed that these phenomenona are quite general, also there were obtained
realizations of such semigroups as semigroups of matrix-valued rational functions
of matrix argument. In all these cases we have double cosets of an infinite-dimensional
analog $G$ of reductive groups with respect to an infinite dimensional analogs $K$ of compact
(generally, non-maximal)
subgroups (as $\U(\infty)\supset\OO(\infty)$). In a recent paper \cite{Ner-finite} on
an infinite-dimensional Grassmannians over finite fields there arises a relatively
new situation. We discuss only a multiplication of double cosets which arises in this context.

\sm

{\bf \punct Result of the paper.} Now we realize the group $G=\GL(\infty,\k)$ as 
a group of matrices infinite to left, right, up, and down. We represent its elements as block matrices
of size $\infty+m+\infty$. Denote by $\P\subset\GL(\infty,\k)$ the group of matrices
of form
$$
\begin{pmatrix}\alpha&\phi&\theta\\
0& \1_m&\psi\\
0&0&\beta
\end{pmatrix}
.$$

We define $\star$-multiplication on double cosets $\P\setminus\GL(\infty,\k)/\P$ as follows.
Consider two elements of $\fra,\frp\in \GL(\infty,\k)$.
We split $\infty+m+\infty$ as $(\infty+N)+m+(N+\infty)$.
Take a sufficiently large $N$ such that these elements 
has the following structure
\begin{equation}
\fra=
\begin{pmatrix}
\1_\infty&0&0&0&0\\
0&a&b&c&0\\
0&d&e&f&0\\
0&g&h&k&0\\
0&0&0&0&\1_\infty
\end{pmatrix}
,\quad
\frp=
\begin{pmatrix}
\1_\infty&0&0&0&0\\
0&p&q&r&0\\
0&v&v&w&0\\
0&x&y&z&0\\
0&0&0&0&\1_\infty
\end{pmatrix}
\label{eq:dve}
.
\end{equation}
Then their $\star$-product is defined by
{\scriptsize%\footnotesize
\begin{multline}
\fra \star \frp=\\
\begin{pmatrix}
\1&0& 0&|&0&|&0&0&0\\
0&\1&0&|&0&|&0&0&0\\
0&0&a&|&b&|&c&0&0\\
\text{---}&\text{---}&\text{---}&&\text{---}&&\text{---}&\text{---}&\text{---}\\
0&0&d&|&e&|&f&0&0\\
\text{---}&\text{---}&\text{---}&&\text{---}&&\text{---}&\text{---}&\text{---}\\
0&0&g&|&h&|&k&0&0\\
0&0&0&|&0&|&0&\1&0\\
0&0&0&|&0&|&0&0&\1
\end{pmatrix}
\begin{pmatrix}
\1&0&0&|&0&|&0&0&0\\
0&p&0&|&q&|&0&r&0\\
0&0&\1&|&0&|&0&0&0\\
\text{---}&\text{---}&\text{---}&&\text{---}&&\text{---}&\text{---}&\text{---}\\
0&v&0&|&v&|&0&w&0\\
\text{---}&\text{---}&\text{---}&&\text{---}&&\text{---}&\text{---}&\text{---}\\
0&0&0&|&0&|&\1&0&0\\
0&x&0&|&y&|&0&z&0\\
0&0&0&|&0&|&0&0&\1
\end{pmatrix}
\label{eq:big}
.
\end{multline}
}
These matrices and their product have size
$$
\infty+N+N+m+N+N+\infty=(\infty+N+N)+m+(N+N+\infty)
.
$$
We unite blocks $\infty+N+N=\infty$, $N+N+\infty=\infty$
and get a matrix of size $\infty+m+\infty$.

We also can write our operation in the following form. Split
$\infty+m+\infty$ as $(\infty+k+k)+m+(k+k+\infty)$, the $m\times m$-block
is the same. Denote
{\small
\begin{equation}
J_k=
\begin{pmatrix}
\1_\infty&0&0&|&0&|&0&0&0\\
0&0&\1_k&|&0&|&0&0&0\\
0&\1_k&0&|&0&|&0&0&0\\
\text{---}&\text{---}&\text{---}&&\text{---}&&\text{---}&\text{---}&\text{---}\\
0&0&0&|&\1_m&|&0&0&0\\
\text{---}&\text{---}&\text{---}&&\text{---}&&\text{---}&\text{---}&\text{---}\\
0&0&0&|&0&|&0&\1_k&0\\
0&0&0&|&0&|&\1_k&0&0\\
0&0&0&|&0&|&0&0&\1_\infty
\end{pmatrix}
\label{eq:Jm}
.
\end{equation}
}
Then 
$$
\fra\star \frp=\fra J_N\frp J_N
.
$$
Notice that the formula 
\begin{equation}
\fra\,\wt\star\, \frp=\fra  J_N \frp 
\label{eq:wt-circ}
\end{equation}
determines the same double coset. 

\begin{theorem}
\label{th}
 The $\star$-multiplication is a well-defined  associative operation
  on $\P\setminus \GL(\infty,\k)/\P$. 
\end{theorem}

Associativity is obvious. Another statement is proved in the next section.

\sm

We define an involution 
$\fra\mapsto\fra^*$ on $\P\setminus \GL(\infty,\k)/\P$ by $\fra\mapsto\fra^{-1}$,
where $\fra^{-1}$ denotes the inverse matrix. Clearly, the operation is well defined

\begin{proposition}
\label{pr}
For any $\fra$, $\frp\in \P\setminus \GL(\infty,\k)/\P$, we have 
$$
(\fra\star \frp)^*=\frp^*\star \fra^*
$$
\end{proposition}

{\sc Proof.} This follows from formula (\ref{eq:wt-circ}).
\hfill $\square$

\sm

Notice that $\fra \star \fra^*=\1$ only for $\fra$ having the form
$$
\fra=\begin{pmatrix}
\1_\infty&0&0\\
0&e&0\\
0&0&\1_\infty
\end{pmatrix}
$$
(the size of the matrix is $\infty+m+\infty$).

\sm

{\sc Remark.}
Formally, the construction is defined for any field (and even for rings).
However, representations of
$\GL(\infty,\k)$ with $\P$-fixed vectors  arise naturally in \cite{Ner-finite}
for $\k$ being finite field.

\section{Proof of Theorem \ref{th}}

\COUNTERS

%{\bf\punct Proof of Proposition \ref{pr}.}

It is sufficient to prove the following statement

\begin{lemma}
Let $\fra$, $\frp\in \GL(\infty,\k)$, $\Phi\in \P$. Then

\sm

{\rm a)} There exists $\Gamma\in \P$ such that 
 $(\fra\cdot \Phi)\star \frp=(\fra\star \frp)\cdot\Gamma$.
 
\sm

{\rm b)} There exists $\Delta\in \P$ such that 
 $\fra\star(\Phi\circ \frp)=\Delta\cdot(\fra\star \frp)$.
\end{lemma}

By Proposition \ref{pr}, it is
sufficient to prove the first statement.

Notice that first and last columns and the first and last rows in matrices
in formulas (\ref{eq:dve}), (\ref{eq:big}), (\ref{eq:Jm}) keep no information.
To reduce size of formulas we omit these columns and rows from the matrices.
In this notation, we have 
$$
\fra\star\frp=
\begin{pmatrix}
p&0&q&0&r\\
bu&a&bv&c&bw\\
eu&d&ev&f&ew\\
hu&g&hv&j&hw\\
x&0&y&0&z
\end{pmatrix}
.
$$

Next, it is sufficient to prove the lemma for $\Phi$ ranging
in a collection of generators
$$
\begin{pmatrix}
\alpha&0&0\\
0&\1_m&0\\
0&0&\beta
\end{pmatrix}
,\quad
\begin{pmatrix}
\1_N&\phi&0\\
0&\1_m&0\\
0&0&\1_N
\end{pmatrix}
,\quad
\begin{pmatrix}
\1_N&0&\theta\\
0&\1_m&0\\
0&0&\1_N
\end{pmatrix}
,\quad
\begin{pmatrix}
\1_N&0&0\\
0&\1_m&\psi\\
0&0&\1_N
\end{pmatrix}
$$
of the group $\P$. Here the size of matrices is $N+m+N$, the matrices
$\alpha$, $\beta$ range in $\GL(N,\k)$; matrices $\phi$, $\psi$, $\theta$
are arbitrary matrices (of appropriate size).

 We examine these generators case by case. First,
$$
\biggl[\fra \cdot \begin{pmatrix}
\alpha&0&0\\
0&\1&0\\
0&0&\beta
\end{pmatrix} \biggr] \star\frp=
\begin{pmatrix}
p&0&q&0&r\\
bu&a\alpha&bv&c\beta&bw\\
eu&d\alpha&ev&f\beta&ew\\
hu&g\alpha&hv&j\beta&hw\\
x&0&y&0&z
\end{pmatrix}
= (\fra\star \frp) 
\begin{pmatrix}
\1&0&0&0&0\\
0&\alpha&0&0&0\\
0&0&\1&0&0\\
0&0&0&\beta&0\\
0&0&0&0&\1\\
\end{pmatrix}
$$
Second, 
\begin{multline*}
\biggl[\fra\cdot \begin{pmatrix}
\1&\phi&0\\
0&\1&0\\
0&0&\1
\end{pmatrix} \biggr]\star\frp
=
\begin{pmatrix}
p&0&q&0&r\\
bu+a\phi u &a&b+a\phi v&c&bw+a\phi w\\
eu+d\phi u&d &e+d\phi v&f&ew+d\phi w\\
hu+g\phi u&g &h+g\phi v&j&hw+g\phi w\\
x&0&y&0&z
\end{pmatrix}
=\\=
(\fra\star \frp) \cdot
\begin{pmatrix}
\1&0&0&0&0\\
\phi u&\1&\phi v&0&\phi w\\
0&0&\1&0&0\\
0&0&0&\1&0\\
0&0&0&0&\1\\
\end{pmatrix}
,\end{multline*}
the $(N+N)+m+(N+N)$-matrix in the right-hand side is contained in the subgroup $\P$.
Next,
\begin{multline*}
\biggl[\fra\cdot \begin{pmatrix}
\1&0&\theta\\
0&\1&0\\
0&0&\1
\end{pmatrix} \biggr]\star\frp
=
\begin{pmatrix}
p&0&q&0&r\\
bu&a&bv&c+a\theta &bw\\
eu&d&ev&f+d\theta &ew\\
hu&g&hv&j+g\theta &hw\\
x&0&y&0&z
\end{pmatrix}
=\\=
(\fra\star\frp)
\cdot
\begin{pmatrix}
\1&0&0&0&0\\
0&\1&0&\theta&0\\
0&0&\1&0&0\\
0&0&0&\1&0\\
0&0&0&0&\1\\
\end{pmatrix}
\end{multline*}
For the last generator, we have
\begin{equation}
\biggl[\fra\cdot \begin{pmatrix}
\1&0&0\\
0&\1&\psi\\
0&0&\1
\end{pmatrix} \biggr]\star\frp
=
\begin{pmatrix}
p&0&q&0&r\\
bu&a&bv&c+b\psi &bw\\
eu&d&ev&f+e\psi &ew\\
hu&g&hv&j+h\psi  &hw\\
x&0&y&0&z
\end{pmatrix}
\label{eq:psi}
\end{equation}
Denote
$$
\begin{pmatrix}
p&q&q\\
u&v&w\\
x&y&z
\end{pmatrix}^{-1}
=:
\begin{pmatrix}
P&Q&R\\
U&V&W\\
X&Y&Z
\end{pmatrix}
$$
Then the right-hand side of (\ref{eq:psi}) is
$$
\begin{pmatrix}
p&0&q&0&r\\
bu&a&bv&c&bw\\
eu&d&ev&f&ew\\
hu&g&hv&j&hw\\
x&0&y&0&z
\end{pmatrix}\cdot
\begin{pmatrix}
\1&0&0& Q\psi&0\\
0&\1&0&0&0\\
0&0&\1&V\psi&0\\
0&0&0&\1&0\\
0&0&0&Y\psi&\1\\
\end{pmatrix}
$$
(the first factor is ). To verify this, we must evaluate 4th column
of the product. We get
\begin{align*}
&pQ\psi+qV\psi+rY\psi=(pQ+qV+rY)\psi=0\\
&buQ\psi+bvV\psi+ c+bwY\psi=c+b(uQ+vV+wY)\psi=c+b\psi\\
&euQ\psi+evV\psi+ f+ewY\psi=f+e(uQ+vV+wY)\psi=f+e\psi\\
&huQ\psi+hvV\psi+ j+hwY\psi=j+h(uQ+vV+wY)\psi=j+h\psi\\
&xQ\psi+yV\psi+zY\psi=(xQ+yV+zY)\psi=0,
\end{align*}
and this completes the proof.

{\tt Math.Dept., University of Vienna,

Oskar-Morgenstern-Platz 1, 1090 Wien,
 Austria

\&

Institute for Theoretical and Experimental Physics,

Bolshaya Cheremushkinskaya, 25, Moscow 117259,
Russia

\&

Mech.Math.Dept., Moscow State University,

Vorob'evy Gory, Moscow

e-mail: neretin(at) mccme.ru

URL:www.mat.univie.ac.at/$\sim$neretin

wwwth.itep.ru/$\sim$neretin
}

\end{document}